\input amstex

\magnification=1200

\input amstex

\centerline{\bf Combinatorial Yamabe Flow on Surfaces}

\bigskip

\centerline{\bf Feng Luo}
\bigskip

\centerline{\bf abstract}

\bigskip

In this paper we develop an approach to conformal geometry of piecewise flat metrics on manifolds. In particular, we formulate the combinatorial Yamabe problem for piecewise flat metrics. In the case of surfaces, we define the
combinatorial Yamabe flow on the space of all piecewise flat metrics associated to a triangulated surface. We show that the flow either develops removable singularities or converges exponentially fast to a constant combinatorial curvature metric. If the singularity develops, we show that the singularity is always removable by a surgery procedure on the triangulation.
  We conjecture that after finitely many such surgery 
changes on the triangulation, the flow converges to the constant combinatorial curvature metric as time approaches infinity.

\bigskip

\S 1. {\bf Introduction}

\bigskip
 
1.1. In
this  paper we consider the class of piecewise flat metrics on a triangulated manifold and develop an approach to
piecewise linear (PL) conformal geometry of these metrics. We also formulate the combinatorial Yamabe problem for piecewise flat metrics.
The notion of conformal class of a metric was originated from Riemannian metrics.
Given a Riemannian metric $g_{ij}$ on a smooth closed manifold, the conformal  change of the metric is a new Riemannian metric of the form $ug_{ij}$ where $u$ is some smooth positive function defined on the manifold. The Yamabe problem states that 
there exists a constant scalar curvature metric in the conformal class of any 
Riemannian metric. This was solved affirmatively for closed Riemannian manifolds in [Ya], [Tu], [Au] and [Sc]. The solution was considered as a milestone in application of non-linear partial differential equations in geometry. 
In trying to develop the analogous PL conformal geometry, one faces the task of defining the notion of
PL conformal class of a PL metric. The most naive approach is to take a local (infinitesimal) point of view. Namely, since the conformality in
smooth category means infinitesimal invariance of angles, thus  PL conformal geometry should preserve the 
of measurment of angles in the PL metrics even at the singularities. However, we are not able to come up with
any reasonable results from this point of view. It is probably due to the that fact that local geometry of smooth
metrics is governed by linear algebra and the local geometry of PL metric is governed by the geometry of simplexes.
We propose to approach the analogous PL conformal geometry of PL metrics from global point of view in this paper. 
Namely, the conformal class of a smooth metric comes from the action of smooth functions on smooth metrics. Thus
the PL conformal class of a PL metric should be a class of PL metrics related by the action of functions defined
on the set of all vertices. Our main results indicate that the PL conformal geometry captures some of the main
features of the smooth conformal geometry on surfaces. 
The situation in PL theory is more complicated in the sense that there are
triangulations on surfaces so that any PL metric associated to the triangulatinos has
none constant PL scalar curvature. This is mainly caused by
the fact that a PL metric is supported by many different triangulations.  
In the case of surfaces, we prove that the existence of a constant PL scalar curvature metric
is a purely topological condition on the triangulation (theorem 1.1). We propose to find the constant PL scalar curvature metric in the PL conformal class by a system of ordinary differential equations. These equations seems to be the right analogy of the Yamabe flow. We establish some of the basic properties of the flow. Especially, we prove
that the flow will not develop essential singularities in finite time and converges exponentially fast to the constant
PL curvature metric if no singularity develops.  The most interesting property of the flow is that when the flow develops singularity, the flow naturally suggests a way to do surgery on the underlying triangulation.
This indicates that the flow tends to search for both the best piecewise flat metric and
the underlying triangulation.

\bigskip

1.2. We begin with a dictionary between piecewise linear (PL) and smooth theory. 
In the smooth theory, we start with a smooth manifold $N^n$ and whereas in the PL theory, we start with a triangulated manifold $(M^n, T)$ where $T$ is the triangulation. The natural analogy of the smooth functions on $N^n$ are the functions defined on the set of all vertices $V$ in the triangulation $T$. Similarly, the analogy of
tensors on $N^n$ are functions defined on the set of all $i$-simplexes in the triangulation. For instance, the de Rham cohomology and simplicial cohomology are
such examples. We define a \it PL metric \rm associated to the triangulation to be a positive real valued 
function defined on the set $E$ of all edges (i.e., 1-simplexes) of $T$ so that
for each n-simplex in $T$, the restriction of the function to the 1-simplexes of
the n-simplex is the edge lengths of some Euclidean  n-simplex. The last condition
is the realizability condition which is given by the Cayley-Menger matrices for
Euclidean n-simplexes.
We may think of a PL metric as a metric on the manifold so that its restriction to
each simplex is isometric to a Euclidean simplex. 
Just like the space of all Riemannian metrics on a smooth manifold is convex, it
can be shown that the space of all PL metrics associated to a triangulation is
a convex set. Indeed, it is an observation of Igor Rivin that the set
$\{ d^2: E \to \bold R |$ where $d: E \to \bold R_{>0}$ is a PL metric associated to $T$\}
is convex in $\bold R^E$.  For any Riemannian metric, the curvature tensor
 assigns to every four tangent vectors at a point a number. For the PL metric, there are
many notions of curvatures (see [CMS],  [CE] and others). The most
natural curvature seems to be the function $K$ which assigns to each (n-2)-simplex 
$\sigma_{n-2}$ the number
$ K(\sigma_{n-2}) = 2\pi - \alpha$ where $\alpha$ is the sum of all dihedral angles of $n$-simplexes
at the (n-2)-simplex. It is shown in [St] that  if the $K(\sigma_{n-2}) >0$ for  all (n-2)-simplexes, then the manifold has finite fundamental group. This is
the analogy of the Myers' theorem in the smooth theory. The \it combinatorial scalar curvature \rm (or PL scalar curvature) of a PL metric $d$ is a function $S_d: V \to \bold R$ whose value at the vertex $v$ is 
$S_d(v) = \sum_{ \sigma_{n-2} >v} K(\sigma_{n-2}) vol(\sigma_{n-2})$ where the summation is over all $(n-2)$-simplexes having $v$ as
a vertex. If $n=2$, then, as a convention, we set $vol(v) =1$.  Note that if $\lambda$ is a positive number
then $S_{\lambda d} = \lambda^{n-2} S_d$ and also that 
 the total scalar curvature  $\sum_{v \in T^0} S(v)$ is exactly
Regge's proposed approximation to the Einstein action (see [Re]). Furthermore, the PL scalar curvature
$S_d$ depends not only on the metric $d$ (as a metric in point set topology) but also on the choice of the underlying triangulation $T$.
If the manifold is a surface, then the scalar curvature $S(v)$ is the curvature $K(v)$.

\bigskip

Using this dictionary, it is now natural to define a \it PL conformal factor \rm for a PL metric 
to be a positive function $u$ defined on the set $V$ of all vertices. Given a PL metric $d: E \to \bold R_{>0}$, we define the PL conformal  change of $d$ to be the new metric
$u*d: E \to \bold R_{>0}$ where $u*d(v v') = u(v)u(v') d(v v')$ and $v v'$ is the
edge with vertices $v$ and $v'$ so that $u*d$ is realizable on each $n$-simplex. The set of
all PL metrics of the form $u*d$ for a fixed $d$ is called the \it PL conformal class \rm of the
PL metric $d$.

\bigskip

The \it combinatorial Yamabe problem \rm asks if there is a constant PL scalar curvature metric within the PL conformal class of
any PL metric. As we will see, there are topological obstructions for the existence of the
constant scalar curvature PL metrics associated to a triangulation even in dimension 2 (see theorem 1.1).
 The goal is to establish the
existence of constant PL scalar curvature metric when the topological obstructions vanish.

\bigskip

1.3. In the case of a triangulated surface, we propose to approach the combinatorial Yamabe problem using a system of ordinary
differential equations. Here is the setup.  Let $T$ be a triangulation of a closed surface $M$ and $V$ and $E$ be the
set of all vertices and edges in $T$. A PL metric associated to $T$ is a positive function $d: E \to \bold R$
so that it satisfies the triangular inequalities on the three edges of any triangle in $T$. Fix the PL metric
$d$ on the surface. Let $u: V \to \bold R_{>0}$ be a conformal factor. For simplicity, let us write
$V =\{ v_1, ..., v_N\}$ and $u_i = u(v_i)$. If $v_iv_j$ is an edge in $T$, then $d(v_iv_j) $ is denoted by
$d_{ij}$.  We define the
\it combinatorial Yamabe flow \rm to be the following system of ordinary differential equations.

\bigskip

$$\frac{du_i}{dt}  = -u_i K_i  \tag 1.1$$
$$ u_i(0)=1$$
Here the curvature $K_i$ is the $K$ curvature of the metric $u*d$ at time $t$ at the i-th vertex $v_i$.
\bigskip

The Gauss-Bonnet theorem for PL metrics states that $\sum_{i=1}^N K_i =2\pi \chi(M)$. Thus the
average curvature $K_{av}$ is $2\pi\chi(M)/N$. In particular, the constant combinatorial curvature PL metric
has curvature $K_{av}$ at all vertices. The combinatorial Yamabe problem asks if this metric exists in the PL conformal class.

There is a combinatorial obstruction for the existence of the constant
curvature PL metric associated to a triangulation $T$. In fact the following holds.
(For a finite set $X$, we use $|X|$ to denote the number of elements in $X$.)

\bigskip

{\bf Theorem 1.1.} \it Fix a triangulation $T$ of a closed topological surface $M$. There exists
a constant PL scalar curvature metric associated to $T$ if and only if for any proper subset
$I$ of the vertices $V$ of $T$, 
$$ |F_I|/|I| > |F|/|V| \tag *$$
where $F$ is the set of all triangles in $T$ and $F_I$ is the set of all triangles
having a vertex in $I$. 

Furthermore, the condition $(*)$ holds for all triangulations of surfaces of non-negative Euler characteristic.\rm

\bigskip

For simplicity, we call a triangulation which supports a constant curvature PL
metric \it admissible. \rm

There are none admissible triangulations on each closed surface of negative Euler characteristic.
It can be shown that given any triangulation of a surface, there is a subdivision of it which
is admissible. However, we are not sure if the subvision can be choosen to be an iterated barycentric subdivision
or even the 2-dimensional subdivision which replace each triangle by four triangles formed by the
barycenters of the vertices and edges.  On the other hand, if our conjecture in section 1.3 holds,
then every triangulation of the surface is equivalent to an admissible one having the same number of
vertices. 

\bigskip

Our main results concerning the combinatorial Yamabe flow can be summarized in the following theorems.

\bigskip

{\bf Theorem 1.2.} \it Fix a piecewise flat metric $d$ on a triangulated surface $(M, T)$. The following holds.

(a). Under the combinatorial Yamabe flow, the curvature evolves according to a combinatorial heat equation of the form

$$ \frac{dK_i}{dt} = \sum_{j=1}^N c_{ij} K_j \tag 1.2$$
where the matrix $[c_{ij}]_{N \times N}$ is symmetric and semi-negative definite. In particular, the total curvature
$\sum_{i=1}^N K_i(t)^2$ is decreasing in time $t$.

(b). The combinatorial Yamabe flow is variational. If we change the variable $u_i$ to $w_i =\log u_i$, then the
combinatorial Yamabe flow is the negative gradient flow of a locally convex function in $w=(w_1, ..., w_N)$.

(c). (local rigidity)  The curvature map $K: \{ u=(u_1, ..., u_N) \in \bold R^N_{>0} | \prod_{i=1}^N u_i =1,$
 $u*d$ is again a PL metric  \} $\to \{ (k_1, ..., k_N) | \sum_{i=1}^N k_i =2\pi \chi(M)\}$ sending a conformal factor $u$ to the curvature
of the metric $u*d$ is a local homeomorphism. 

\rm
\bigskip

We do not know if this local rigidity can be improved to be a global rigidity result.

\bigskip

To understand the long time behavior of the flow, we introduce the normalized equation

$$ \frac{du_i}{dt} = -u_i (K_i -K_{av})  \tag 1.3$$
$$ u_i(0)=1$$

The two equations (1.1) and (1.3) are equivalent in the sense that $u_i(t)$ solves (1.1) 
if and only if $ e^{2\pi \chi(M) t/N } u_i(t)$ solves (1.3). Furthermore, if 
$u(t) =(u_1(t), ..., u_N(t))$ is a solution to the normalized equation (1.3) then the product $\prod_{i=1}^N u_i(t) =1$ for all time $t$.
To understand the asymptotic behavior of the solution of (1.3), we have to analysis the potential formation
of singularities in the normalized equation (1.3).
There are only two types of singularities which may occur. Suppose $u(t)
$ is a solution of (1.3) in the time interval $[0, L)$ where $L \leq \infty$. We say
the solution develops an \it essential singularity \rm at time $L$ if there is an index
$i$ and a sequence of time $t_n$ approaching $L$ so that
$\lim_{t_n \to L} u_i(t_n) = 0$. We say the solution develops \it removable singularity \rm at time $L$ if
there is a sequence of time $t_n$ approaching $L$ so that
$ u_i(t_n)$ remains in a compact set in $\bold R_{>0}$ for all $i$ and there
is a triangle $\Delta v_iv_jv_k$ in $T$ which degenerates into a line segment as $t_n \to L$.
The last condition means that the triangular inequality $u_i d_{ik} u_k + u_kd_{kj} u_j > u_id_{ij} u_j$
becomes equality as $t_n \to L$ for some indices $i,j,k$. 

\bigskip

{\bf Theorem 1.3.} \it 
(a) For any triangulated surface, the normalized combinatorial Yamabe flow will not develop essential singularity
in finite time. 

(a) If the triangulation is admissible, the normalized combinatorial Yamabe flow  will not develop essential singularity
at time infinity. 
\rm
\bigskip

If removable singularity occurs, say the vertex $v_k$ is moving toward the
interior of the edge $v_iv_j$, then we change the original triangulation $T$ as follows. Let
$v_l$ be the fourth vertex in $T$ so that $v_iv_jv_l$ forms a triangle in $T$. Then we replace the
triangulation $T$ by a new triangulation $T'$ which is obtained from $T$ by deleting the edge
$v_iv_j$ and adding a new edge $v_kv_l$. In particular, two triangles $\Delta v_iv_jv_k$ and $\Delta v_iv_jv_l$
in $T$ are replaced by two new triangles $\Delta v_k v_l v_i$ and $\Delta v_k v_l v_j$. This is the most common surgery operation on
triangulations of surfaces. After the combinatorial surgery on $T$, we run the combinatorial Yamabe
flow on $T'$ with initial metric coming from the final stage of metric on $T$ at time $L$.
We conjecture that after finitely many such surgery operations on the triangulation $T$, the combinatorial
Yamabe flow converges exponentially fast to the constant curvature metric. This is supported by the following.

\bigskip

{\bf Theorem 1.4.} \it If no singularity develops in  the normalized combinatorial Yamabe flow, then the solution
converges exponentially fast to a constant curvature PL metric as time approaches infinity. \rm

\bigskip

Under the normalized flow, the distance function restricted to the 1-skeleton of
the triangulation will stay in the same quasi-isometric class in any compact time interval in $\bold R_{>0}$.
The geometry of the quasi-isometry class changes only when the essentialy singulariy develops at time infinity.

If we start with a single Euclidean triangle, and define the curvature $K_i$ at the i-th vertex to be $\pi -\theta_i$
where $\theta_i$ is the inner angle at the i-th vertex, then it can be shown that the corresponding normalized
equation (1.3) has solution for all time and converges exponentially fast to the equilateral triangle as
time tends to infinity.

\bigskip

The following is likely to hold for the combinatorial Yamabe flows on surfaces.

\bigskip

{\bf Conjecture.} \it (a) A constant PL scalar curvature metric is unique in its  PL conformal class.

(b)  The combinatorial Yamabe flow will not develop essential singularities on any triangulation.

(c) The combinatorial Yamabe flow will converge to the constant curvature PL metric after finite number of
surgeries on the triangulation. \rm

\bigskip

There are some evidences indicating that the combinatorial Yamabe flow will develop removable singularities for
some PL metrics on admissible triangulations. But we do not have the explicit examples yet.

It seems there is an interesting similarity between the singularity formation in combinatorial Yamabe flow on surfaces and
Hamilton's Ricci flow program in dimension 3.
 \bigskip

1.4. There are several interesting questions concerning the combinatorial Yamabe problem in higher dimension. To carry out
the same program in higher dimension, the key ingredient that is missing is the following local rigidity
property. 

\bigskip

{\bf Question.} \it  Given an Euclidean n-simplex $\sigma_n$ and an (n-2)-face $\sigma_{n-2}$ of $\sigma_n$, let $a(\sigma_{n-2}, \sigma_n)$ be the dihedral angle
of the n-simplex at the (n-2)-simplex.  Now suppose $d: E \to  \bold R_{>0}$ is the
edge lengths of the n-simplex. Consider conformal factors $u: V \to \bold R_{>0}$ and
the new PL metric $u*d$ on the n-simplex. Define the function 
$S_{u*d}: V \to \bold R^n$ to be $S_{u*d}(v) = \sum_{ \sigma_{n-2} > v} a(\sigma_{n-2}, \sigma_n) vol(\sigma_{n-2})$ in the metric $u*d$.  
 Then for fixed metric $d$, $S_{u*d}$ is
a smooth function of $u$. Is the rank of the Jacobian matrix of the map $S_{u*d}$ considered as a function defined in an open set in $\bold R^V$ to $\bold R^V$ always equal
to $n-1$? \rm

\bigskip

An affirmative answer to this question for $n \geq 3$ will give a strong evidence that the
higher dimensional combinatorial Yamabe flow exists, i.e., one would call equations of the type
$\frac{du_i}{dt} = f(u_i) S(v_i)$ the combinatorial Yamabe flow where $f(x)$ is some universal no-where zero function depending
on the dimension $n$. The function $f$ is so chosen that the corresponding matrix $[\frac{\partial S_i}{\partial u_j }f(u_j)]$
is symmetric.  This question is also related to the following. Given a PL metric
$d$ associated to a triangulation, consider the set of all conformal factors $u$
for $d$. Let $R(u*d) = \sum_{\sigma_{n-2}} K(\sigma_{n-2}) vol(\sigma_{n-2})$ be
the Regge action of the metric $u*d$. Is the Hessian of $R(u*d)$ considered
as a function of $u$ semi-positive definite?

\bigskip

The other related question we are  considering now is the combinatorial Yamabe flow in hyperbolic or spherical background metric, i.e.,
we use hyperbolic or spherical simplexes instead of the Euclidean ones.  The work of [CL] suggests that in general
it is easier to work in hyperbolic geometry.

\bigskip

There are many interesting questions in the field of PL metric theory. For instance what is the
right analog of the Ricci curvature? What is the Laplacian operator on the space of all functions
defined on the $i$-simplexes? 
Is it true that positive PL scalar curvature carries topological information? To be more
precise, suppose $M$ is a 3-manifold which supports a positive PL scalar curvature. Is it ture that the fundamental group of $M$ contains no none-abelian surface group? 
In the smooth case, the result was proved by Schoen and Yau [SY].
 
\bigskip

1.5. 
In [CR], Cooper and Rivin define the concept of combinatorial scalar curvature differently. 
 In their definition, the combinatorial scalar curvature is invariant under the scaling of  metric.  On the other hand, the scalar curvature in Riemannian geometry has the property that $S(kg_{ij}) = 1/k S(g_{ij})$ for any positive constant k. (In our case, the PL scalar curvature satisfies $S(kg) = k^{n-2} S(g)$.) It can be shown every topological 3-manifold has a triangulation and an associated PL metric whose PL scalar curvature is positive in the sense of [CR]. Thus the positivity of PL scalar curvature defined in the sense of [CR] does not carry any topological information. However, at each vertex,
the scalar curvature defined in [CR] is a measurment of the infinitesimal rate of change of the volume of a small ball centered at the vertex. Thus, it is a combinatorial analogy of the scalar curvature in the infinittesimal sense. 
 The work [Ga] follows the approach of [CR] and defines a combinatorial Yamabe flow for ball packing metrics.
The PL scalar curvature defined in this paper is more closely related to Regge's calculus.
\bigskip

1.6.
The paper is organized as follows. In \S2, we establish some basic properties of the geometry of triangles. In \S3, \S4, and S5, we prove 1.1, 1.2, 1.3, and 1.4. In the appendix A,
we provide a detail calculation involved in the proof of theorem 2.1. In appendix B, we
use the feasible flow theorem for network flow to establish the necessity part of theorem 1.1.

\bigskip

1.6. This work is supported in part by the NSF. We thank 
Ben Chow, X-S. Lin and D. Sullivan for discussions. We thank the referee for making some nice suggestions.

\bigskip

\S 2. {\bf Geometry of Euclidean Triangles}

\bigskip

For simplicity, let  $\Delta =\{ (x_1, x_2, x_3) \in \bold R_{>0}^3| x_i + x_j > x_k, $ where i, j, k are
distinct\} be the set of points whose coordinates satisfy the triangular inequalities. We use $\bold R_{>0}$ to denote the set of all positive numbers, and $\Delta v_1 v_2v_3$ to denote the
triangle having vertices $v_1, v_2, $ and $v_3$. If $X$ is a finite set, we use $|X|$ to denote the number of elements in 
$X$.

\bigskip

2.1. Fix three numbers $d_i, d_j ,d_k$ satisfying triangular inequalities, i.e.,
$(d_i, d_j ,d_k) \in \Delta$ and  three vertices $v_i, v_j ,v_k$. Consider the conformal factor $u:
\{v_i, v_j, v_k\} \to \bold R_{>0}$ so that $(d_i/u_i, d_j/u_j, d_k/u_k) \in \Delta$ where $u_r = u(v_r)$. Define $x_r = d_ru_su_t$ where $\{r,s,t\}=\{i,j,k\}$.
Due to the choice of $u_r$'s, $(x_i, x_j, x_k) \in \Delta $.
Construct a Euclidean triangle $\Delta v_iv_jv_k$ so that the length of the edge opposite to $v_r$ is
$x_r$. Let $\theta_r$ be the inner angle at the vertex $v_r$. Then $\theta_r$ is a smooth function
of $(u_i, u_j, u_k)$. 

\bigskip

{\bf Theorem 2.1.} \it The $3 \times 3$ matrix $[ \frac{\partial \theta_r }{\partial u_s} u_s]_{3 \times 3}$ is symmetric, semi-negative
definite and has rank 2 whose null space is $\{ (t,t,t) \in \bold R^3|$ $t \in \bold R\}$.  Furthermore, if we let
$a_{rr} = x_r^2$ and $a_{rs} = -x_r x_s \cos (\theta_t)$ where $\{r,s,t\}=\{i,j,k\}$ and let $A$ be
the area of the triangle $\Delta v_iv_jv_k$, then $[ \frac{\partial \theta_r }{\partial u_s} u_s]_{3 \times 3}$ $=-\frac{1}{2A} [a_{rs}]_{3 \times 3}$.
\rm
\bigskip

The proof is a straight forward computation. We defer it to the appendix A. It is easy to see that
the matrix $[a_{rs}]_{3 \times 3}$ is symmetric and semi-positive definite. Indeed, since $x_i = x_j \cos \theta_k + x_k
\cos \theta_j$, we see that the sum of entries in every row is zero. Thus $det([a_{rs}]_{3 \times 3}) =0. $  On the
other hand, the diagonal entries are positive so are the determinants of the principal $2 \times 2$ submatrices.
Thus the $3\times 3$ matrix $[a_{rs}]_{3 \times 3}$ is semi-positive definite whose null space is the diagonal in
$\bold R^3$.

\bigskip

{\bf Corollary 2.2.} \it The matrix $[\frac{\partial \theta_r }{\partial u_s} u_s]_{3 \times 3}$ depends only on the
three inner angles $\theta_i, \theta_j$ and $\theta_k$ of the triangle. In particular, if the conformal factor
$u =(u_i, u_j, u_k)$ varies in a region so that the inner angles of the triangle $\Delta v_iv_jv_k$ lie in a compact set in the open interval $(0, \pi)$, then there is
a positive constant $\lambda$ so that the negative eigenvalues of the matrix $[\frac{\partial \theta_r}{\partial u_s} u_s]_{3 \times 3}$
are less than $-\lambda$ for all $u$ in the region. \rm

\bigskip

2.2. Using the same notations as above, let $w_r = \log u_r$ for $r=i,j,k$. Then 
the $3 \times 3$ matrix $[ \frac{\partial \theta_r }{\partial w_s}]_{3 \times 3}$ is symmetric and semi-negative definite of
rank 2 whose null space is $\{(t,t,t) \in \bold R^3| t \in \bold R\}$.
Consider the non-convex space $W =\{ w=(w_i, w_j , w_k) \in \bold R^3 | (d_i e^{-w_i}, d_j e^{-w_j}, d_k e^{-w_k}) \in
\Delta \}$. The space $W$ is simply connected since it is the image of the convex
space $\{ (r_i,r_j,r_k) \in \bold R_{>0}^3| (d_ir_i, d_jr_j, d_kr_k ) \in \Delta \}$  under the homeomorphism $h (r_i, r_j, r_k)  = (-\log r_i, -\log r_j,- \log r_k)$.
Furthermore, if $w \in W$, then $w + (t,t,t)$ is again in $W$ for any real number $t$. Take
any vector $(a_i, a_j ,a_k) \in \bold R$   and let $\Omega = (a_i- \theta_i) dw_i + (a_j -\theta_j) dw_j +
(a_k - \theta_k) dw_k$
be a smooth 1-form defined on the space $W$. It is closed due to the symmetry of the matrix
$[ \frac{\partial \theta_r }{\partial w_s}]_{3 \times 3}$. Define a smooth function $F(w): W \to R$ by

$$ F(w) = \int_0 ^w \Omega. \tag 2.1$$

This is well defined since the form $\Omega$ is closed and the space $W$ is simply connected.

\bigskip
 
{\bf Corollary 2.3.} \it The function $F: W \to \bold R$ is a locally convex function so that it becomes
a locally strictly convex function when restricted to the planes $\{ (w_i, w_j, w_k) \in W | w_i + w_j + w_k $ is
a constant \}. Furthermore, 

(a). if $a_i+a_j+a_k =\pi$, then $F(w + (t,t,t)) = F(w)$ for any real number $t$.

(b). if $(a_i,a_j,a_k) \in \bold R^3_{>0}$ and
$a_i+a_j + a_k =\pi$, then for any sequence of points  $w^{(n)}$   in $W$,
$ \limsup_{n \to \infty} F(w^{(n)}) =\infty $ if and only if
 $ \limsup_{n \to \infty}
 \max_{ \{r,s\} \subset \{  i,j,k\}}   (| w^{(n)}_r - w^{(n)}_s|) = \infty$. \rm

\bigskip

{\bf Proof. }  The local convexity is due to the fact the Hessian of $F$ is
the matrix $-[\frac{\partial \theta_r }{\partial w_s}]_{3 \times 3}$ which is semi-positive definite.
Furthermore, since the matrix is positive definite on the linear space $w_i+w_j+w_k=0$, it follows that $F$ is locally strictly convexity on the planes.

To see (a), it suffices to verify that $\int_w ^{w + (d,d,d)} \Omega =0$.
Now take the line segment $t w + (1-t)(d,d,d)$ in $W$ to evaluate the line integral
$\int_w ^{w+(d,d,d)} \Omega$. We find the integral
becomes $\int_0^1 (( a_i+a_j+a_k )-(\theta_i + \theta_j +\theta_k)) dt =0$ since
the integrant is 0.

To see (b), let $\Pi: \bold R^3 \to \{(w_1,w_2, w_3) \in \bold R^3 | w_1+w_2+w_3=0\}$ be the
orthogonal projection. Then due to the condition $ai+a_j+a_k = \pi$ and part (a), we have
$F(w) = F(\Pi(w))$. Also, we have 
$ \limsup_{n \to \infty} \max_{ \{r,s\} \subset \{  i,j,k\}}   (| w^{(n)}_r - w^{(n)}_s|) =\infty $ 
if and only if their projections $\Pi(W^{(n)})$ is unbounded. 
Thus if $ \limsup_{n \to \infty}
\max_{ \{r,s\} \subset \{  i,j,k\}}   (| w^{(n)}_r - w^{(n)}_s|) $ is finite, then the
function $F(w^{n})$ remains bounded.

Now suppose $\limsup_{n \to \infty} \max_{\{r,s\} \subset\{i,j,k\}}|w^{(n)}_r - w^{(n)}_s| = \infty$. By taking
a subsequence if necessary, we may assume without loss of generality that
$\lim_{n \to \infty} (w^{(n)}_i - w^{(n)}_j) = \infty$. By part (a), we may further assume 
after adding $w^{(n)}$ by a diagonal vector $(d_n, d_n, d_n)$ that
$w^{(n)}$ is in the hypersurface $\{ (w_i, w_j, w_k) \in W |  e^{-w_i} + e^{-w_j}
+  e^{-w_k} = 1\}$. In particular, $w^{(n)}_r \geq 0$. Thus, $\lim_{n \to \infty}
w^{(n)}_i = \infty$. 
By 
 the constraints $e^{-w_i} + e^{-w_j} + e^{-w_k} =1$
and the triangular inequalities $ d_i e^{-w_i} + d_r e^{-w_r} > d_s e^{ -w_s}$
for $\{r,s\} =\{j,k\}$, both $w^{(n)}_j$ and $w^{(n)}_k$ are bounded from above.
 This shows that the triangle $\Delta v_iv_jv_k$ with edge lengths $d_i e^{w_j^{(n)} + w_k^{(n)}},   d_j e^{ w_i^{(n)} + w_k^{(n)}}$ and $ d_k e^{ w_i^{(n)}+ w_j^{(n)}}$
degenerates into a half-line since two of the edge lengths tend to infinity and the third
remains bounded.
Thus $\theta_i(w^{(n)})$ tends to zero. Therefore the dominate term in 
the integration $\int_0^{w^{(n)}} \Omega$ is $ \int_0^{w^{(n)}_i} (a_i -\theta_i) dw_i$ since the other two integrals are bounded. Because $a_i>0$ and $\lim_n \theta_i(w^{(n)}) =0$, we see that $\lim_{n \to \infty} F(w^{(n)}) =
\infty.$  QED

\bigskip

2.3. {\bf Remarks.} 
The situation in PL conformal geometry is very similar to the approach to circle packing by Colin de Verdiere in [Cv]. The complication in the combinatorial Yamabe problem is caused by the fact that the space $W$ is not convex and the function $F: W \to \bold R$ is not proper.

\bigskip

\S3. {\bf A Proof of Theorem 1.1}

\bigskip

That the condition $(*)$ is sufficient was established in [CL] theorem 1.1 where we
take the weight function $\phi =0$. Indeed, condition (1.3) for $\phi=0$ in theorem 1.1 in [CL] is exactly the condition $(*)$.
Thus in this case, there is a very special type PL metric associated
to the triangulation whose curvature at each vertex is $2\pi\chi(M)/|V|$. The metric
is obtained by assigning to each vertex $v_i$ a positive number $r_i$ and define the length at the edge $v_i v_j$ to be $r_i + r_j$.

To show that the condition is also necessary, we use the feasible flow theorem for
network flow. Since the method of the proof is irrelevant to the rest of the paper, 
we defer the proof to appendix B. 

The fact that the condition $(*)$ holds for all triangulations of surfaces of non-negative Euler characteristic was established in [CL].

\bigskip

\S4. {\bf Proofs of Theorems 1.2. and 1.3}

\bigskip

4.1.
To derive the evolution of the curvature $K_i$ for the combinatorial Yamabe flow $\frac{du_i}{dt}  = -u_i K_i$,
we note that the evolution of an individual
inner angle, say $\theta_i^{jk}$ in a triangle $\Delta v_iv_jv_k$ at the vertex $v_i$, is the following.

$$\frac{d \theta_i^{jk}}{dt}  = \sum_{r=i,j,k}  \frac{\partial \theta_i^{jk} }{\partial u_r} \frac{du_r}{dt}  $$
$$ = -\sum_{r=i,j,k} \frac{\partial \theta_i^{jk}}{\partial u_r} u_r K_r \tag 4.1$$

Now the curvature $K_i$ at the i-th vertex is $2\pi -\sum_{j,k} \theta_i^{jk}$ where the summation is over all
triangle $\Delta v_iv_jv_k$ in $T$.
 Thus
$\frac{ d K_i}{dt} = - \sum_{j,k} \frac{d \theta_i^{jk}}{dt} $.  By (4.1) and theorem 2.1, we see that
the curvature evolution equation is

$$\frac{ dK_i}{dt}  = \sum_{r=1}^N   \sum_{j,k}^N  \frac{\partial \theta_i^{jk}}{\partial u_r } u_r  K_r = \sum_{r=1}^N c_{ir} K_r.$$
We can write the coefficient matrix as a sum of matrices, 
$$  [c_{ir}]_{N \times N} = \sum_{j,k} [\frac{\partial \theta_i^{jk}}{\partial u_r} u_r]_{N \times N}$$
where the sum is over all ordered edges $\{j,k\}$.  The matrix
$A^{j,k} =  [\frac{\partial \theta_i^{jk}}{\partial u_r} u_r]_{N \times N}$ 
 is symmetric and semi-positive
definite with only one $3 \times 3$ submatrix being non-zero and all other entries being zero by   theorem 2.1.
Thus the matrix $[c_{ij}]$ is symmetric and semi-positive definite. Furthermore, by the same identity above and theorem 2.1, we see that
the rank of $[c_{ij}]$ is $N-1$ whose null space is spanned by the vector $[1, 1, ..., 1]$.

Take the time derivative of $G(t) =\sum_{i=1}^N K_i^2(t)$. We get
$$ \frac{dG}{dt}  = 2 \sum_{i=1} K_i K_i' = 2 \sum_{i, j=1}^N c_{ij}K_iK_j \leq 0.$$
Thus the total curvature square is decreasing along the solution. This establishes
part (a) of theorem 1.2.

\bigskip

4.2.
To show part (b), if we change the variable $u_i$ to $w_i = \log u_i$, then the
combinatorial Yamabe flow becomes $\frac{dw_i}{dt}  = - K_i$. Note that the symmetry of the matrix
$[\frac{\partial \theta_r}{\partial u_s} u_s]_{3 \times 3}$ implies that the $N \times N$ matrix
$[\frac{\partial K_r}{\partial w_s}]_{N \times N}$ is symmetric and is semi-positive definite.
Its rank is (N-1) and its null space is $\{(t,t,...,t) | t \in \bold R\}$.

 Define  $\Cal W$ to be the space
$\{ w=(w_1, ..., w_N) \in \bold R^N | $ if $v_iv_jv_k$ forms a triangle in $T$, then
$(w_i, w_j, w_k) \in W$\}. Note that $\Cal W$ is not a convex space but is simply
connected since it is the image of the convex set $\{ (r_1,..., r_N) \in \bold R^N_{>0}
|$ if $v_iv_jv_k$ forms a triangle in $T$, then $d_ir_i + d_jr_j > d_kr_k$ \} under the
diffeomorphism $(r_1, ..., r_N) \to (-\log r_1, ..., -\log r_N)$. Also, if
$w \in \Cal W$, then $w+(t,t,..., t)$ is still in $\Cal W$.

Define a 1-form $\Omega = \sum_{i=1}^N K_i dw_i$ on the space $\Cal W$. This
form is closed since the matrix $[ \frac{\partial K_r}{\partial w_s}]_{N \times N}$ is symmetric. Thus the function
$$ F(w) = \int_0^w \Omega$$
is well defined on the simply connected space $\Cal W$. The function $F$ is locally convex since its Hessian is the semi-positive
definite matrix $[\frac{\partial K_r}{\partial w_s}]_{N \times N}$ and is locally strictly convex on the planes
$\Cal W \cap \{ w | w_1 + ... + w_N =$ constant \} since the matrix is positive definite when restricted to the
plane $\{w | w_1 + ...+ w_N =0\}$. Due to $\frac{\partial F}{w_i} = K_i$, it follows 
 that the combinatorial
Yamabe flow $\frac{dw_i}{dt}  = -K_i $ is the negative gradient flow of
a locally convex function defined on $\Cal W$.

\bigskip

4.3.
To see the local rigidity, we note that under the transformation from $u$ to $w$, the
space $u$-space becomes the (N-1)-dimensional smooth manifold $P =\Cal W \cap \{ w_1+... +w_N =0\}$. The curvature map $K$ becomes the map sending $w \in P$ to the gradient
of a locally strictly convex function $F|_P$ defined on $P$. Due to the local strictly convexity,
the gradient map is locally injective. Since the image of $K$ is 
$\{(k_1,..., k_N) | \sum_{i=1}^N k_i = 2\pi \chi(M)\}$ which is also a manifold of
dimension (N-1), by the invariance of domain theorem, the curvature map $K$ must
be a local homeomorphism. This establishes part (c). 

\bigskip

4.4. {\bf Remark}. It can be shown that the curvature map $K$ is a local diffeomorphism.
However, we are not able to show that is globally a diffeomorphism.

\bigskip

4.5. To prove  part (a) of theorem 1.3, note that 
if $u_i(t)$ solves the combinatorial Yamabe flow, then due to the universal
curvature bound that $  (2-|E|) \pi \leq K_i < 2\pi$, we have
 $$ 1/c e^{-ct} \leq u_i(t) \leq c e^{ ct}$$
for some positive constant $c$ on the time interval $[0, L)$ where the solution exists.
Thus no essential singularity develops on any triangulation in finite time.

\bigskip

4.6. To prove part (b) of theorem 1.3, since the triangulation is admissible, we find a 
PL flat metric associated to the triangulation so that its curvature is constant
$2\pi \chi(M)/N$. Let $a_i^{jk}$ be the inner angle at vertex $v_i$ in triangle
$\Delta v_iv_jv_k$ in the metric.  
For the normalized equation $u_i' = -u_i(K_i -K_{av})$, we make a change of
variable $w_i = \log u_i$ as before. By the normalization condition that
$\prod_{i=1} u_i (t) =1$ for all time, we have $\sum_{i=1}^n w_i(t) = 0$.
By the same calculation, we still have the fact that the $N \times N$ matrix
$[ \frac{\partial (K_r-K_{av})}{\partial w_s}]_{N \times N}$ is again symmetric,
semi-positive definite of rank (N-1) whose null space is $\{(t,t,...,t) \in
\bold R^N | t \in \bold R\}$.

Now consider the space $\Cal W $ as in subsection 4.2.
We construct a specific 1-form $\Omega$ as
$$ \Omega = \sum_{i=1}^N (\sum_{j,k} ( a_i^{jk} -\theta_i^{jk}) dw_i ) \tag 4.2$$
where $\theta_i^{jk}$ is the inner angle of the metric $u*d$ at the vertex
$v_i$ inside the triangle $\Delta v_iv_jv_k$. By the choice of the inner
angles $a_{i}^{jk}$, the 1-form $\Omega$ is exactly $\sum_{i=1}^N (K_i -K_{av}) dw_i$.

Let us now consider the associated function $F(w) = \int_0^w \Omega$ defined
on $\Cal W$. This function is well defined since $\Cal W$ is simply connected and
$\Omega$ is closed. Furthermore, by the construction, the combinatorial Yamabe flow
$w_i'= -(K_i -K_{av})$ is the negative gradient flow of $F(w)$. In particular,
if $w(t)$ is a solution to the combinatorial Yamabe flow, the function $h(t) = F(w(t))$ is
decreasing in time $t$.

On the other and, we can rewrite the summation in (4.2) as follows, 

$$ F(w) = \sum_{\Delta v_iv_jv_k \in T}   \int_0^{w} (a_i^{jk} -\theta_i^{jk}) dw_i + 
(a_j^{ik} -\theta_j^{ik})dw_j + (a_k^{ij} -\theta_k^{ij})dw_k \tag 4.3$$ where the sum is over all triangles in
$T$.

Since $a_i^{jk} + a_j^{ki} + a_{k}^{ij} = \pi$ over each
triangle, by corollary 2.3 (a), $F(w + (d,..., d)) = F(w)$ for all $w \in
\Cal W$. Now suppose the normalized combinatorial Yamabe flow $w_i'  = -(K_i -K_{av})$
develops an essential singularity at time $t=L \leq \infty$. Then, due to the normalization equation $\sum_{i=1}^N w_i(t) =0$, there is a sequence of
time $t_n \to L$ and two indices $i,j$ so that
$\lim_{t_n \to \infty} w_i(t_n) = \infty$ and
$\lim_{t_n \to \infty} w_j(t_n) = -\infty$.  

We claim in this case $\lim_{t_n \to L} F(w(t_n)) = \infty$. This will contradict
the fact that $F(w(t))$ is decreasing in time $t$.

To see that claim, we use the same argument used in the proof of corollary 2.3(b). 
 Indeed,
after adding $w(t_n)$ by a vector of the form $(d_n,...d_n)$, we may assume that
$w(t_n)$ is in the subspace $\{ (w_1, ..., w_N) \in \bold R^N | $ 
$e^{-w_1} + e^{-w_2} + ...+ e^{-w_N} = 1$\}. In this case, we have
$w_j(t_n) \geq 0$ for all $j$. Also there are two indices $i,j$ so that

$$\lim_{t_n \to L} w_i(t_n) = \infty  \quad \text{and} \quad \quad w_j(t_n) \quad \text{ remains bounded.} \tag 4.4$$

Now by (4.3) the integral $F(w(t_n))$ is the sum of finitely many integrals of type
(2.1) over each triangle $\Delta v_iv_jv_k$ where the vector $(a_i^{jk},
a_j^{ik}, a_k^{ij})$ satisfies the condition in corollary 2.3(b).
Thus, by the same argument used in the proof of corollary 2.3(b), the
integral $ \int_0 ^{w(t_n)} (a_i^{jk} -\theta_i^{jk}) dw_i + 
(a_j^{ik} -\theta_j^{ik})dw_j + (a_k^{ij} -\theta_k^{ij})dw_k$ is either bounded
or tends to infinity depending on $\limsup_{n \to \infty} \max_{\{r,s\} \subset \{i,j,k\}} ( | w_r(t_n) - w_s(t_n)|)$
is finite or infinite. However, the infinite case must occur due to (4.4). This shows that
$F(w(t_n))$ tends to infinity.
QED

\bigskip

\S5. {\bf A Proof of theorem 1.4.}

\bigskip
By the assumption, the solution $u(t) = (u_1(t), ..., u_N(t))$ of the normalized combinatorial Yamabe flow
exists for all time so that there are no singularities forming at time equal to
infinity. This means that $u_i(t)$'s are in some compact interval
in $\bold R_{>0}$ and also all inner angles $\theta_i^{ij}(t)$ are in some compact
interval inside the interval $(0, \pi)$. By corollary 2.2,  this implies that the there is
a positive constant $\lambda$ so that the eigenvalues of 
coefficient matrix $[c_{rs}]_{N \times N}$ considered as a bilinear form restricted to the
subspace $\{w \in \bold R^N | w_1+ ...+ w_N = 0\}$ is always bounded by $-\lambda$ for all time $t \in [0, \infty)$, i.e.,

$$ \sum_{i,j} c_{ij} w_i w_j \leq -\lambda \sum_{i} w_i^2, \quad  \text{when $\sum_{i=1}^N w_i = 0$.} \tag 5.1$$

To prove the theorem, it suffices to show that the curvature $K_i(t)-K_{av}$ converges to 0 exponentially fast, i.e., there is a positive constants $c_1, c_2$ so
that 

$$ | K_i(t) -K_{av}| \leq c_1 e^{ -c_2t}.  \tag 5.2$$

Assuming (5.2) holds, then we can solve $u_i(t)$ from the combinatorial Yamabe flow 
and get $$ u_i(t) = e^{ - \int_0^t (K_i(s) -K_{av}) ds}.$$  This shows that
$\lim_{t \to \infty} u_i(t)$ is a positive real number for all indices. Thus
the metrics $u(t) d$ converges to the constant PL curvature metric.

To establish (5.1), let us consider the function
$$ G(t) = \sum_{i=1}^N (K_i(t) -K_{av})^2$$

Its derivative can be calculated as

$$G'(t) = 2 \sum_{i,j} c_{ij} (K_i -K_{av}) (K_j -K_{av}) \tag 5.3$$

By (5.1), we have $G'(t) \leq -\lambda G(t)$. Thus $G(t) \leq C e^{-\lambda t}$. This establishes
(5.2). QED

\bigskip

5.2. {\bf Remark.} One can now show that the normalized combinatorial Yamabe flow for a
single Euclidean triangle $\frac{du_i}{dt}  = -(\pi/3 - \theta_i) u_i$ converges to the
equilateral triangle as follows. By the same argument as in the proof of theorem 1.3,
we see that essential singularity never occur. If a removable singularity occurs
at time $t=L$, then the three inner angles tend to $0, 0, \pi$. On the other hand the
sum of the square of the curvature $\sum_{r} (K_i-K_{av})^2$ achieves its supremum 
value only when the inner angles are $0,0,\pi$. Thus by theorem 1.2(a) that the
sum of the square of the curvature is decreasing, removable singularities never occur. 
By theorem 1.4, we see the solution converges. This ends the proof.
\bigskip

\centerline{\bf Reference}

\bigskip

[A] Aubin, T.,  
Équations différentielles non linéaires et problème de Yamabe concernant la courbure scalaire. 
J. Math. Pures Appl. (9) 55 (1976), no. 3, 269--296.

[CE] Cao, J.-G.; Escobar, J. F., A New 3-dimensional Curvature Integral Formula for PL-manifolds of Non-positive Curvature, preprint, Feb. 2000.

[CL] Chow, B., and Luo, F., Combinatorial Ricci flows on surfaces, preprint, 2002, 

http://front.math.ucdavis.edu/math.DG/0211256.

[CMS] Cheeger, J.; M\"uller, W.; Schrader, R., On the curvature of piecewise flat spaces. Comm. Math. Phys. 92 (1984), no. 3, 405--454.  

[CR]  Cooper, D.; Rivin, I., Combinatorial scalar curvature and rigidity of ball packings. Math. Res. Lett. 3 (1996), no. 1, 51--60. 

[Cv] Colin de Verdière, Y., Un principe variationnel pour les empilements de cercles.  Invent. Math. 104 (1991), no. 3, 655--669.

[FF] Ford, L. R., Jr.; Fulkerson, D. R., Flows in networks. Princeton University Press, Princeton, N.J. 1962

[Ga] Glickenstein, D., A maximum principle for combinatorial Yamabe flow, preprint, 2002.

[Re] Regge, T., General relativity without coordinates, Nuovo Cimento, 19,(1961),558-571.

[Sc] Schoen, R., Conformal deformation of a Riemannian metric to constant scalar curvature. J. Differential Geom. 20 (1984), no. 2, 479--495. 

[St] Stone, D. A.,
Geodesics in piecewise linear manifolds. 
Trans. Amer. Math. Soc. 215 (1976), 1--44.

[SY] Schoen, R.; Yau, S. T., 
Existence of incompressible minimal surfaces and the topology of three-dimensional manifolds with nonnegative scalar curvature. 
Ann. of Math. (2) 110 (1979), no. 1, 127--142.

[Tu] Trudinger, N. S., 
Remarks concerning the conformal deformation of Riemannian structures on compact manifolds. 
Ann. Scuola Norm. Sup. Pisa (3) 22 1968 265--274.

[Ya] Yamabe, H., 
On a deformation of Riemannian structures on compact manifolds. 
Osaka Math. J. 12 1960 21--37.

\bigskip

\centerline{\bf Appendix A. A Proof of Theorem 2.1}

\bigskip

We will carry out the computational aspect of the proof of theorem 2.1.

\bigskip

The following lemma was established in [CL], Lemma A-1.

\bigskip

\noindent
{\bf Lemma A-1.} \it Suppose $\Delta v_iv_jv_k$ is a Euclidean triangle of area $A$ so that
the inner angle at $v_i$ is $\theta_i$ and the length of the edge $v_jv_k$ is $x_i$. Then
$\theta_i$ is a smooth function of $(x_i, x_j, x_k)$.  The partial derivatives of the function
are given by, 

(a) $\frac{\partial  \theta_i}{\partial x_i} = x_i/(2A) $.

(b) $\frac{\partial \theta_i}{\partial x_j} = - \frac{\partial \theta_i}{\partial x_i} \cos(\theta_k)$.

\rm

\bigskip

In our case, we fix a set of positive numbers $d_i, d_j, d_k$ and choose a conformal factor
$(u_i, u_j, u_k) \in \bold R_{>0}^3$. The edge lengths of the triangle are $x_i = d_i u_ju_k$. Thus $\frac{\partial x_i }{\partial u_i} =0$
and $\frac{\partial x_i}{\partial u_j} = x_i/u_j$.

Now the partial derivative can be calculated by the chain rule,

$$\frac{\partial \theta_i }{\partial u_r} = \sum_{s} \frac{\partial \theta_i }{\partial x_s} \frac{ \partial x_s }{\partial u_r}$$

$$= \sum_{ s \neq r} \frac{\partial \theta_i}{\partial x_s } \frac{x_s}{u_r}$$

This shows that $\frac{\partial \theta_i}{\partial u_r} u_r = \sum_{s \neq r}\frac{ \partial \theta_i }{\partial x_s }x_s$. Now
use the lemma A-1 above together with the fact that $x_i = x_j \cos \theta_k + x_k \cos \theta_j$, we
obtain $\frac{\partial \theta_i}{\partial u_r} u_r = -\frac{a_{ir}}{2A} $ as stated in theorem 2.1.

\bigskip

\centerline{ \bf Appendix B. An Application of the Feasible Flow Theorem to Theorem 1.1 }

\bigskip

We now verify that the condition $(*)$ in theorem 1.1 is also sufficient for the existence of constant curvature PL metric
associated to the triangulation. This proof follows the same ideas appeared in [Cv].

Let us begin with the feasible flow theorem for network flow. Suppose $G =(V, E)$ is a directed graph, i.e., a graph so that each edge is oriented. Here $V$ is the set of all vertices and $E$ is the set of all oriented edges. For each oriented edge $x \in E$, let $int(x)$ and $end(x)$ be the initial vertex and the end vertex
of the edge $x$. (It is possible that $int(x) = end (x)$.) For any subset $I \subset V$, let $int(I)
=\{ x \in E | int(x) \in I,  end(x) \notin I\}$, and let $end(I) =\{ x \in E | end(x) \in I, int(x) \notin I\}$.

Assume there is a lower capacity bound $a: E \to [-\infty, \infty]$ and an upper capacity bound $b: E \to [-\infty, \infty]$ so that
$a(x) \leq b(x)$ for all $x \in E$.

A feasible flow on the graph $G$ is a function $\phi: E \to (-\infty, \infty)$ so that Kirchoff's current law is
satisfied, i.e., for each vertex $v$

$$ \sum_{ x \in end(v)} \phi(x) = \sum_{x \in int(v)} \phi(x) $$
and $a(x) \leq  \phi(x) \leq b(x)$ for all $x \in E$.

\bigskip

{\bf Feasible Flow Theorem.} \it A feasible flow exists if and only if for every non-empty subset $U \subset V$
so that $U \neq V$,

$$ \sum_{ x \in end(U)} b(x) \geq \sum_{ x \in int(U)} a(x). \tag B1$$

\rm

We now apply the theorem to show that condition $(*)$ in theorem 1.1 is a necessary condition. To this end,  let
the set of all vertices in the triangulation be $V$ and the set of all triangle in $T$ be $F$. If $v \in V$
and $f \in F$, then $f>v$ means $v$ is a vertex of $f$. Let $z$ be an extra vertex. Define an oriented graph
$G$ as follows. The set of all vertices in $G$, denoted by $G^0$ is $V \cup F \cup \{z\}$. The set of
all oriented edges $G^1$ is $\{ (f,v) | f \in F, v \in V, f>v\} \cup \{ (z,f)| f \in F\} \cup \{ (v,z) | v \in V\}$.

Now suppose there is a constant PL flat metric associated to the triangulation. Define a feasible flow $\phi: G^1 \to
[0, \infty)$
on the graph $G$ as follows: $\phi((f,v))$ is the inner angle of the triangle $f$ at the vertex $v$, and
$\phi((z,f))= \pi$ and $\phi((v,z)) = \pi |F|/|V|$. 
Define the lower capacity $a$ for $G$ as follows:  $ a((f,v)) = \epsilon > 0$ which is smaller than any of the
inner angles, $a(x) =\phi(x)$ for all other edges.
Define the upper capacity $b$ for $G$ as follows:   $b((f,v)) = \infty$, $ b(x) = a(x)$ for all other edges.

To verify the  Kirchoff's law for $\phi$, it suffices to check the following three statements.
At each vertex $v \in V$, it states that the sum of the inner angles at the vertex is  $2\pi - 2\pi (\chi(M)/|V|) = \pi |F|/|V|$, i.e.,
the curvature is constant. At each vertex $f \in F$, it states that the sum of the inner angles of a triangle is $\pi$,
and at the vertex $z$, it is the trivial statement that $ \pi |F| = |V|(\pi |F|/|V|)$.

Recall that for 
a non-empty subset $I$ of $V$ so that $I \neq V$, we use $F_I =\{ f \in F | f > v $ for some $v \in I$\}.
For this set $I$, we consider the proper subset $U = I \cup F_I$ of vertices $G^0$. We claim that the feasibility condition (B1) for this $U$ is
exactly the condition $(*)$ in theorem 1.1.

Indeed, we have $end(U) =\{ (f, v) | v \in I, f \notin F_I,   f>v\}  \cup \{ (z,f) | f \in F_I\} 
$. Note that by the choice of $F_I$,
the set $\{(f,v) |v \in I, f \notin F_I, f>v\} =\emptyset$. Thus $end(U) =\{(z,f) | f \in F_I\}$.
Also, we have
$int(U) =\{ (v, z) | v \in I\} \cup \{ (f,v) | f \in F_I, v \notin I, f>v\}$. Thus the feasibility condition states that

$$ \sum_{f \in F_I} \pi  \geq \sum_{v \in I} \pi (|F|/|V|) + \sum_{\{(f,v)| f \in F_I, v \notin I, f > v\}} \epsilon
> \sum_{v \in I} \pi (|F|/|V|) .$$

Thus we have $|F_I|/|I| > |F|/|V|$.

\bigskip

It can be shown easily that the condition $(*)$ in theorem 1.1 is exactly the feasibility condition in the
feasible flow theorem. Thus one can in fact avoid using the result from [CL] in the proof
of theorem 1.3. Namely, we choose the set of inner angles $a_i^{jk}$ in the proof by applying the feasible flow
theorem on the same graph $G$ and the same upper and lower capacities $a$ and $b$.

\bigskip

Department of Mathematics

Rutgers University

New Brunswick, NJ 08845

email: fluo\@math.rutgers.edu

\end
\end